\documentclass[12pt]{amsart}

\usepackage{amsmath,amssymb,amscd,amsfonts}

\newtheorem{thm}{Theorem}[section]
\newtheorem{lem}[thm]{Lemma}
\newtheorem{rem}[thm]{Remark}
\newtheorem{prop}[thm]{Proposition}

\newtheorem{assu-nota}[thm]{Assumption--Notation}

\theoremstyle{remark}

\newcommand{\C}{\mathbb C}
\newcommand{\Z}{\mathbb Z}

\newcommand{\F}{\mathbb F}
\newcommand{\pp}{\mathbb P}

\DeclareMathOperator{\Pic}{Pic}

\DeclareMathOperator{\Num}{Num}

\DeclareMathOperator{\Hom}{Hom}

\newcommand{\la}{\lambda}
\newcommand{\Ga}{\Gamma}

\newcommand{\Si}{\Sigma}
\newcommand{\si}{\sigma}
\newcommand{\tsi}{\tilde{\sigma}}
\newcommand{\fie}{\varphi}
 \newcommand{\tK}{\tilde{K}}
 \newcommand{\tA}{\tilde{A}}

\newcommand{\calO}{\mathcal{O}}

\newcommand{\inv}{^{-1}}
\newcommand{\lra}{\longrightarrow}

\numberwithin{equation}{section}

\begin{document}
\title{Enriques surfaces with eight nodes}

\author{Margarida Mendes Lopes, Rita Pardini}
\date{}
\begin{abstract}
 A nodal Enriques surface can have at most 8 nodes. We give an explicit
description of Enriques surfaces with 8 nodes, showing that they are quotients of
products of elliptic curves by a group isomorphic to
$\Z_2^2$ or to $\Z_2^3$ acting freely in codimension 1. We use this result to
show that if $S$ is a minimal  surface of general type with $p_g=0$ such that the
image of the bicanonical map is birational to an Enriques surface then $K^2_S=3$
and the bicanonical map is a morphism of degree 2.
\end{abstract}

\maketitle
\section{Introduction} It is well known that a nodal
 Enriques surface has at most 8 nodes. In this note, applying the technique for
the study of nodal surfaces developed in
\cite{nodi},   we are able to characterize completely the Enriques surfaces with
8 nodes.
 We show that every such surface is a  quotient of a product of elliptic curves
by a group isomorphic to
$\Z_2^2$ or to $\Z_2^3$ acting freely in codimension 1 (see Theorem \ref{main}
for the precise statement).

In the last section we apply this classification result to show that if $S$ is a
minimal  surface of general type with $p_g=0$ such that the image of the
bicanonical map $\fie$ is birational to an Enriques surface then $K^2_S=3$ and $\fie$ is
a morphism of degree 2. This result refines Theorem 3 of
\cite{xiaocan}, ruling out one of the possibilities presented there.

The paper is organized as follows: in section 2 we recall the facts we need from
\cite{nodi}; in section 3 we describe in detail the construction of Enriques
surfaces with 8 nodes as quotients of products of elliptic curves; in section 4
we show that the 8 nodes on a nodal Enriques surface form an even set and we
prove the classification theorem; finally, in section 5 we apply the previous
result to the study of the bicanonical image of a surface with $p_g=0$.

\medskip
\noindent{\bf Acknowledgements.} We wish to thank JongHae Keum who kindly
communicated to us a different proof of Lemma 4.2.

\smallskip The present collaboration takes place in the framework of the european
contract EAGER, no. HPRN-CT-2000-00099. The first author is a member of CMAF and
of the Departamento de Matem\'atica da Faculdade de Ci\^encias da Universidade de
Lisboa and the second author is a member of GNSAGA of CNR.

\medskip
\noindent{\bf Notation and conventions.} We work over the complex numbers; all
varieties are assumed to be compact and algebraic. We do not distinguish between
line bundles and divisors on a smooth variety, using the additive and the
multiplicative notation interchangeably. Linear equivalence is denoted by
$\equiv$ and numerical equivalence by $\sim_{num}$. $\Num(X)$ is the group of line
bundles of the variety  $X$ modulo numerical equivalence.
As it is usual, we denote by
$\kappa(X)$ the Kodaira dimension of a variety $X$ and by $\rho(X)$ the Picard
number of  $X$.

We recall that a {\em  nodal Enriques surface} $\Sigma$ is a normal projective
regular surface  whose singular  points are nodes and such that
$K_{\Sigma}\not\equiv 0$, $2K_{\Sigma}\equiv 0$.

\section{Nodal surfaces and codes}\label{nodal} In this section we recall the
basic facts that we will need  about nodal surfaces and we establish the notation.

A {\em  nodal surface} $\Sigma$ is a normal projective surface whose singular
points are {\em nodes}, i.e. they are  singularities analytically isomorphic to
the hypersurface singularity $x^2+y^2+z^2=0$.  In particular,
$\Si$ has canonical singularities and $K_{\Si}$ is  Cartier.  Denote by
$P_1\dots P_k$ the nodes of
$\Sigma$ and let
$\eta\colon Y\to
\Si$ be the minimal resolution. One has $K_Y=\eta^*K_{\Si}$,
$\chi(\calO_Y)=\chi(\calO_\Si)$ and for every
$i=1\dots k$ the curve $C_i:=\eta\inv P_i$ is a {\em nodal curve}, namely a
smooth rational curve such that $C_i^2=-2$, $K_YC_i=0$.

Conversely, given a surface $Y$ and a disjoint set $C_1\dots C_k$ of nodal curves
of $Y$, there exist a nodal surface $\Sigma$ and a birational morphism
$\eta\colon Y\to \Si\ $ that contracts $C_1\dots C_k$ to nodes $P_1\dots P_k$ and
is an  isomorphism on the complement of $\cup _iC_i$. Thus it is equivalent to
consider  the nodal surface $\Si$ or the smooth surface $Y$ together with the set
of disjoint nodal curves $C_i$.

The geometry of a nodal surface is often studied by means of the corresponding
{\em  binary code}.
 Recall that a binary code  of length $k$ is a  linear subspace  of
$\F_2^k$.  Given a vector
$v=(x^1\dots x^k)$  of the code,
 the \emph{weight} of $v$ is the number of indices $i$ such that $x^i\ne 0$. The
code $V$ associated to the set of disjoint nodal curves $C_1\dots C_k\subset Y$
is the kernel of the   homomorphism
$\psi\colon \F_2^k\to \Pic(Y)/2\Pic(Y)$ defined by $(x^1\dots x^k)\mapsto \sum
x^i [C_i]$, where $[D]$ denotes the  class of a divisor $D$. In other words,
$v=(x^1\dots x^k)$ is in
$V$ if and only if $\sum_ix^iC_i$ is divisible by $2$ in $\Pic(Y)$, namely if and
only if there is a line bundle
$L$ on $Y$ such that $2L\equiv \sum_ix^iC_i$. If for a nonempty subset $I$ of
$\{1\dots k\}$ the element  $\sum_{i\in I}C_i$  is divisible by $2$ in $\Pic(Y)$,
we say that $\{P_i\ |i\in I\}$ is an {\em even} set of nodes. Hence the nonzero
elements of $V$ correspond to the  even sets of nodes of $\Si$.

Replacing
$\Pic(Y)$ by
$\Num(Y)$ and linear equivalence by numerical equivalence, one defines in the same
way the code
$V_{num}$. The nonzero vectors of $V_{num}$ correspond to the {\em numerically
even} sets of nodes of $\Si$.  Clearly,
$V$ is contained in
$V_{num}$ and the two codes may or may not be equal. The code $V_{num}$ can also
be described in the following way. Denote by
$\Ga'\subset \Num(Y)$ the lattice spanned by $C_1\dots C_k$ and by $\Ga$ its
primitive closure, namely the smallest primitive sublattice of $\Num(Y)$
containing $\Ga'$. The lattice $\Ga$ has rank $k$ and the natural
map  $V_{num}\to \Ga/\Ga'$ defined by $(x^1\dots x^k)\mapsto
\frac{1}{2}\sum_ix^iC_i$ is an isomorphism. If we denote by
$\Delta$ the discriminant of
$\Ga$ then we have the following relation:
\begin{equation}\label{discriminant} 2^k=2^{2\dim V_{num}}\Delta.
\end{equation}

Let $v=(x^1\dots x^k)\in V_{num}$  and write $\sum _ix^iC_i\sim_{num}2L$ for a
suitable line bundle
$L$ of $Y$.  Since
 $K_Y L=0$,  $L^2$ is even by the adjunction formula and thus the weight of
$v$ is divisible by $4$. We say that a curve $C_i$ {\em appears in} a subspace
$W$ of $V$ (or $V_{num}$) if $W$
 is not contained in the subspace
$\{x^i=0\}$ of $\F_2^k$.

The following theorem from \cite{nodi} shows how one can obtain information on
the geometry of a nodal surface from the code associated to the set of nodes.

We recall that a Galois cover is said to be {\em totally ramified} if the Galois
group  is generated by the elements with nonempty fixed locus.
\begin{thm}\label{codecover}
 Let $\Si$ be a nodal surface and let $V$ be  the corresponding code. Let
$W\subset V$ be a subcode of dimension $r$ and let
$m$ be the number of nodes of $\Si$ that appear in $W$. Then there exists a
totally ramified Galois cover $\pi\colon Z\to \Si$ such that:
\begin{itemize}
\item[i)] the Galois group of $\pi$ is $G:=\Hom(W,\C^*)$;
\item[ii)] $\pi$ is branched precisely on the nodes of $\Si$ that appear in $W$;
\item[iii)] $Z$ is a nodal surface and the singular set of $Z$ is the inverse
image of the nodes of $\Si$ that do not appear in $W$;
\item[iv)] the invariants of $Z$ are the following:
$$\kappa(Z)=\kappa(\Si);\quad
\chi(\calO_Z)=2^r\chi(\calO_\Si)-m2^{r-3};\quad
K^2_Z=2^rK^2_{\Si}.$$
\end{itemize}
\end{thm}
\begin{proof} See \cite{nodi}, Proposition 2.1 and Proposition 2.3. Both
propositions are stated only  for  the case
$W=V$ but the proofs extend verbatim to the general case.
\end{proof}
\section{The examples}\label{examples}

Recall that a nodal surface $\Si$ is an Enriques surface if the minimal
desingularization
$Y$ of $\Si$ is an Enriques surface. It is well known that a nodal Enriques
surface has at most 8 nodes (see, e.g., \cite{mi}). The following are two
examples of Enriques surfaces with 8 nodes.
\smallskip

\noindent{\bf Example 1.} Let $D_1$, $D_2$ be elliptic curves and let
$a\in D_1$,
$b\in D_2$ be  nonzero points of order 2. Denote by $e_1$, $e_2$ the standard
generators of $\Z_2^2$. We let
$\Z_2^2$ act on $D_1$ as follows:
$$x_1\overset{e_1}{\lra} -x_1;\quad x_1\overset{e_2}{\lra} x_1+a.$$

Analogously we let $\Z_2^2$ act on $D_2$ as:
$$x_2\overset{e_1}{\lra} x_2+b;\quad x_2\overset{e_2}{\lra} -x_2.$$ We consider
the diagonal action of $\Z_2^2$ on $A:=D_1\times D_2$,  we set $\Si:=A/\Z_2^2$
and we denote by $\pi\colon A\to \Si$ the quotient map. The singularities of $\Si$
are $8$ nodes that are the images of the 16 fixed points of $e_1+e_2$. The map
$\pi$ is branched precisely on the nodes of $\Si$, hence we have
$\pi^*K_{\Si}=K_A=0$. Thus
 $\Si$ is minimal of zero Kodaira dimension. Considering the action of $\Z_2^2$
on $H^0(A,\Omega^i_A)$,
$i=1,2$ one checks that $p_g(\Si)=q(\Si)=0$ and thus $\Si$ is an Enriques
surface. A similar argument shows that for $i=1,2$ the surface $Z_i:=A/e_i$ is a
smooth minimal bielliptic surface, while
$K:=A/(e_1+e_2)$ is the Kummer surface of $A$.  The projections $A\to D_i$,
$i=1,2$, descend to pencils of elliptic curves $p_i\colon \Si\to
D_i/\Z_2^2=\pp^1$. The singular fibres of $p_i$ are the following: two smooth
double fibres, occurring at the images in $\pp^1$ of the points of order 2 of
$D_i$,  and two singular double fibres, occurring at the images in $\pp^1$ of the
fixed points of
$e_1+e_2$. Each singular double fibre contains 4 nodes of $\Si$ and is supported
on a smooth rational curve, hence it corresponds to a fibre
 of  type
$I_0^*$ on the resolution $Y$ of $\Si$. If we denote by
$f_i$ the class of  a fibre of $p_i$, $i=1,2$, then $f_i$ is divisible by 2 in
$\Pic(\Si)$ and we have $f_1f_2=4$.  The  4 nodes lying on the same fibre of
$p_1$ or
$p_2$ are an even set.   A singular double fibre of  $p_1$ and a singular double
fibre of $p_2$ intersect at 2 nodes. So the 8 nodes of $\Si$ are divided into 4
pairs  such that the union of any 2 such pairs is even.
\medskip

\noindent{\bf Example 2.} Let $D_1$, $D_2$ be elliptic curves and let
$a_i\in D_1$,
$b_i\in D_2$,
$i=1,2,3$, be the  nonzero points of order 2. Denote by $e_1$, $e_2$,
$e_3$ the standard generators of
$\Z_2^3$. We  let
$\Z_2^3$ act on $D_1$ as follows:
$$x_1\overset{e_1}{\lra} x_1+a_1;\quad x_1\overset{e_2}{\lra} x_1+a_2;
\quad x_1\overset{e_3}{\lra} -x_1.$$

Analogously we let $\Z_2^3$ act on $D_2$ as:
$$x_2\overset{e_1}{\lra} x_2+b_1;\quad x_2\overset{e_2}{\lra} -x_2; \quad
x_2\overset{e_3}{\lra} x_2+b_3.$$ We consider the diagonal action of
$\Z_2^3$ on
$A:=D_1\times D_2$, we set $\Si:=A/\Z_2^3$ and we denote by $\pi\colon A\to
\Si$ the quotient map. Arguing as in Example 1, one shows that $\Si$ is a minimal
Enriques surface with $8$ nodes.  The subgroups $G_2:=<e_1, e_2>$ and $G_3:=
<e_1, e_3>$ act freely on
$A$ and the corresponding quotients   are minimal bielliptic surfaces. The
elements
$e_2+e_3$ and $e_1+e_2+e_3$ have 16 fixed points each and the corresponding
quotients are Kummer surfaces.

The projections $A\to D_i$, $i=1,2$, descend to pencils of elliptic curves
$p_i\colon
\Si\to\pp^1$. Both pencils have two smooth double fibres and two singular
double fibres containing 4 nodes each (corresponding to fibres of type $I_0^*\,$
on the resolution
$Y$ of $\Si$). If we denote by
$f_i$ the class of  a fibre of $p_i$, $i=1,2$, then $f_i$ is divisible by 2 in
$\Pic(\Si)$ and we have $f_1f_2=8$. A singular double fibre of $p_1$ and a
singular double fibre of
$p_2$ either intersect in 4 nodes or they meet at 2 smooth points of $\Si$.
\medskip

We now establish some facts that we will need in the following sections.
\begin{lem}\label{pg} Let $\Si$ be a nodal Enriques surface and let $\pi\colon
Z\to\Si$ be a totally ramified Galois cover such that $Z$ has canonical
singularities and such that the branch locus of $\pi$ is contained in the
singular set of $\Si$.

Then
$p_g(Z)=0$.
\end{lem}
\begin{proof} We have $K_Z=\pi^*K_{\Si}$, since $\pi$ is \'etale in codimension
1. Hence $2K_Z\equiv 0$, and $p_g(Z)\ne 0$ iff $K_Z\equiv 0$. Let $K\to \Si$ be
the Kummer cover of $\Si$ and consider the following diagram, obtained by base
change:
\[
\begin{CD} Z' @>>>Z\\ @VVV @VVV\\ K@>>>\Si.
\end{CD}
\] If $K_Z\equiv 0$, then the double cover $Z'$ is the disjoint union of two
connected components, each mapping isomorphically to $Z$. By the commutativity of
the diagram, this shows that $\pi\colon Z\to \Si$ factors through $K\to\Si$,
contradicting the assumption that $\pi$ is totally ramified.
\end{proof}
\begin{prop} \label{lattice} Let $\Si$ be a nodal Enriques surface as in Example
1 or in Example 2 and let $Y\to\Si$ be the minimal desingularization. Let
$f_1=2A_1$, $f_2=2A_2\in
\Pic(Y)$ be the classes of the fibres of the elliptic pencils of $Y$ induced by
$p_1$ and $p_2$, and  let
$\Ga$ be the smallest primitive sublattice of
$\Num(Y)$  containing the classes  of the exceptional  curves of
$Y\to\Si$. Then: $$\Ga^{\perp}=< A_1,A_2>.$$
\end{prop}
\begin{proof} The inclusion $<A_1,A_2>\subset \Ga^{\perp}$ is obvious. Since both
sublattices have rank 2, to prove equality it is enough to show that they have
the same discriminant. If $\Si$ is as in Example 1, then $<A_1,A_2>$ is
unimodular and the result is immediate.

Assume now that $\Si$ is as in Example 2. In this case, the discriminant of
$<A_1,A_2>$ is equal to $4$.  The  discriminant of
$\Ga^{\perp}$ is equal to the discriminant $\Delta$ of $\Ga$, since $\Ga$ is
primitive and the intersection form on $\Num(Y)$ is unimodular. In addition, we
have
$\Delta= 2^{8-2\dim V_{num}}$ by (\ref{discriminant}). Hence we have to show that
$\dim V_{num}\le 3$. Since $V$ is a subspace of $V_{num}$ of codimension at most
1,  we are  going  to show that
$\dim V=2$.

We have seen above that $V$ contains two disjoint even sets $J_1$ and $J_2$ of
order 4, each contained in a singular double fibre of $p_1$ and in a singular
double fibre of
$p_2$. Assume by contradiction that there is an even set $J$  of 4 nodes different
from
$J_1$ and
$J_2$ and denote by $J'$ the complement of $J$. Then each singular double fibre of
$p_1$ and
$p_2$ contains   2 nodes of
$J$ and 2 nodes of $J'$. Let $\pi\colon Z\to \Si$ be  a $\Z_2^2-$cover associated
to the span $W$ of $J$ and $J'$ in $V$ (cf. Theorem
\ref{codecover}).  The surface
$Z$ is smooth  with $\chi(\calO_Z)=0$ by Theorem
\ref{codecover} and it has
$p_g(Z)=0$ by Lemma \ref{pg}. Since $K_Z=\pi^*K_{\Si}$, it follows that
$Z$ is bielliptic and
$2K_Z\equiv 0$. The cover
$\pi$ factors as
$Z\to  Z_1\to \Si$, where
$Z_1\to \Si$ is a double cover  branched over $J$. The pull-back to $Z_1$ of a
singular double fibre of $p_1$ or $p_2$ is again a fibre of the same type, hence
it is not the double of a Cartier divisor of $Z_1$. This shows that the pull-back
on
$Z_1$ of the general fibre of both $p_1$ and $p_2$ is connected. Hence the
pull-backs to $Z$ of the general fibres of
$p_1$ and $p_2$   either are  connected or they are the disjoint union of two
smooth
 elliptic curves. It follows that the intersection number of the two elliptic
pencils of
$Z$ is $\ge 8$, but since $2K_Z\equiv 0$  this contradicts the classification of
bielliptic  surfaces by Bagnera and De Franchis (see e.g. \cite{beauville},
Ch. VI). So in this case $V$ is generated by $J_1$ and $J_2$, and so $\dim
V=2$.
\end{proof}

\section{The classification theorem}\label{classification}

This section is devoted to the proof of our main result:

\begin{thm}\label{main} Let $\Si$ be an Enriques surface with $8$ nodes.\newline
Then there  exist elliptic curves $D_1$, $D_2$ such that $\Si$ is constructed
from $D_1$,  $D_2$ as in Example 1 or in Example 2.\end{thm} We need the
following Lemma.

\begin{lem}\label{even} If $\Si$ is an Enriques surface with 8 nodes, then the
nodes of
$\Si$ are an even set.
\end{lem}
\begin{proof}
 Let $K\to \Si$ be the K3 cover of $\Si$. The surface
$K$ has  16 nodes, hence $K$ is  the Kummer surface of an abelian surface $A$
(see \cite{nikulin}, Theorem 1). We wish to show that the involution
$\si\colon K\to K$ associated to the double cover
$K\to \Si$ can be lifted to an involution $\si'$ of $A$.  Denote by
$\tilde{K}$ the minimal resolutions of the singularities of $K$ and let $C_1\dots
C_{16}$ be the nodal curves of $\tilde{K}$ arising in the resolution of the nodes
of $K$. Taking base change, one gets the following diagram:
\[
\begin{CD} \tilde{A} @>>>A\\ @VVV @VVV\\ \tilde{K}@>>>K
\end{CD}
\] where $\tilde{A}$ is the blow up of $A$ at the 16 points of order 2. Denote by
$\tilde{\si}$ the involution of $\tilde{K}$ induced by $\si$. The map
$\tilde{A}\to \tilde{K}$ is flat of degree 2, branched on $C_1+\dots+ C_{16}$,
and thus there exist a line  bundle $L$ on
$\tilde{K}$ such that $2L\equiv \sum C_i$. Denote by $V(L)$ the total space   of
$L$ and by $p\colon V(L)\to \tK$ the projection. The surface $\tA$ is isomorphic
to the zero locus on  $V(L)$ of the section  $z^2-p^*f$ of $p^*L^2$, where $z$ is
the tautological section of
$p^*L$ and $f$ is a section of $L^2$ vanishing on $\sum C_i$.

Notice that $L$ is determined uniquely by the condition $2L\equiv \sum C_i$, since
$\tilde{K}$ is simply connected. The divisor $\sum C_i$ is preserved by
$\tilde{\si}$, hence $\tsi^*L\cong L$,
$\tsi$ can be lifted to an automorphism $\si_L$ of $V(L)$ and $f$ is an
eigenvector for the action of
$\tsi$ on $H^0(\tilde{A},L^2)$. It follows that, up to composing with an
automorphism  of $L$ lifting the identity of $\tK$, we can assume that $\si_L$
maps $\tA$ to itself. Thus $\tsi$ can be lifted to an automorphism
$\si'$ of
$\tilde{A}$. In turn, $\si'$   induces an automorphism of $A$, that we denote
again by
$\si'$,  that lifts
$\si\colon K\to K$. Notice that $\si'$ acts freely on $A$, since $\si$ acts
freely on
$K$.  The order of
$\si'$ is either $2$ or $4$. Assume that it is $4$, so that
${\si'}^2=-1_A$. If we write
$\si'z=gz+a$,  with $g$  an automorphism of $A$ and $a\in A$, this gives
$g^2=-1_A$. It follows that the eigenvalues of the differential of $g$ at $0$ are
equal to $i$ or to $-i$. Thus the morphism $g-1_A\colon A\to A$ is surjective,
hence there  exists $z_0\in A$ such that
$(g-1_A)z_0=-a$. This is the same as saying that $z_0$ is a fixed point of
$\si'$. Thus we have a contradiction, showing that  ${\si'}^2=1$.

Set $Z:=A/\si'$. Then we have a commutative diagram:
\[
\begin{CD} A @>>>Z\\ @VVV @VVV\\ K@>>>\Si.
\end{CD}
\] The surface $Z$ is smooth, since it is a free quotient of a smooth surface,
and, by the commutativity of the diagram,
$Z\to \Si$ is a double cover  of $\Si$ branched precisely over the 8 nodes of
$\Si$. Thus the nodes of $\Si$ are an even set.
\end{proof}

\begin{proof}[Proof of Theorem \ref{main}]
 By Lemma \ref{even}, we know there exists a double cover  $Z\to
\Si$  branched over the nodes of $\Si$. By Proposition
\ref{codecover},
$Z$ is  smooth,    $\chi(\calO_Z)=0$ and $\kappa(Z)=0$. We
have
$p_g(Z)=0$ by Lemma
\ref{pg}.  The map $\pi$ is \'etale in codimension 1, hence
$K_Z=\pi^*K_{\Si}$ and $2K_Z\equiv 0$. So $Z$ is minimal bielliptic. By the
classification of bielliptic surfaces of Bagnera--De Franchis (cf.
\cite{beauville}, Ch. VI) there exist elliptic curves
$D_1$,
$D_2$ such that
$Z$ is one of the following:
\begin{itemize}
\item[a)] the quotient of $D_1\times D_2$ by the diagonal action of
$\Z_2$, where
$\Z_2$ acts on
$D_1$ by $x_1{\lra} -x_1$,  and on $D_2$ by $x_2{\lra} x_2+b$, with $b$ a point
of order 2 of $D_2$.

\item[b)] the quotient of $D_1\times D_2$ by the diagonal action of
$\Z_2^2$, where the standard generators $e_1$, $e_2$ of
$\Z_2^2$ act on $D_1$ by:
$$x_1\overset{e_1}{\lra} x_1+a_1;\quad x_1\overset{e_2}{\lra} x_1+a_2,$$ where
$a_1\ne a_2$  are points of order 2 of $D_1$, and they  act on $D_2$ by:
$$x_2\overset{e_1}{\lra} x_2+b_1;\quad x_2\overset{e_2}{\lra} -x_2,$$ where
$b_1$ is a point of order 2 of $D_2$.
\end{itemize}

  We consider the composite map
$D_1\times D_2\to Z\to \Si$ and we wish to show that it is a Galois cover with
Galois group isomorphic to
$\Z_2^2$ in case a) and to
$\Z_2^3$ in case b).   As in the proof of Lemma \ref{even}, we have a commutative
diagram:
\[
\begin{CD} A @>>>Z\\ @VVV @VVV\\ K@>>>\Si,
\end{CD}
\] where $K$ is the K3 cover of $\Si$ and $A$ is an abelian surface. The cover
$A\to Z$ is the \'etale  cover given by the canonical class of $Z$, therefore we
have
$A=D_1\times D_2$ in case a) and $A=(D_1\times D_2)/e_1$ in case b).
 Consider case a) first. Here $K$ is the Kummer surface of
$D_1\times D_2$ and the involution  associated with $A\to K$ can be written as
$x_1\to -x_1+a$, $x_2\to -x_2+b'$. In addition, by  a suitable choice of the
origin of
$D_2$, we may assume $b'=0$. This involution and the involution associated with
$D_1\times D_2$  are in the Galois group of $D_1\times D_2\to \Si$, which is
therefore isomorphic to $\Z_2^2$. Using this remark it is easy to see that $a$ is
a point of order
$2$ of $D_1$. In addition,
$a$ must be nonzero, since otherwise the branch locus of $D_1\times D_2\to
\Si$ would have dimension 1. Thus in this case $\Si$ is the surface of Example 1.

Consider now case b).
 Arguing as before, one shows that the Galois group of $D_1\times D_2$ contains
elements of the following form:
$$(x_1,x_2)\overset{e_1}{\lra}(x_1+a_1, x_2+b_1)$$
$$(x_1,x_2)\overset{e_2}{\lra}(x_1+a_2, -x_2)$$
$$(x_1,x_2)\overset{e_3}{\lra}(-x_1, x_2+b_3),$$ with $a_1\ne a_2\in D_1$ and
$b_1\ne b_3\in D_2$ points of order 2. So in this case
$\Si$ is the surface of Example 2.

\end{proof}

\section{The bicanonical image of a surface with $p_g=0$}\label{application} In
this section we apply the previous results to prove the following:
\begin{thm} Let $S$ be a minimal surface of general type with $p_g(S)=0$,
$K^2_S\ge 3$ and let
$\fie\colon S\to X\subset \pp^{K^2_S}$ be the bicanonical map  of $S$. If $\fie$
is not birational, then either

i) $X$ is a rational surface,

 or

ii) $K^2_S=3$,
$\fie$ is a morphism of degree 2 and $X\subset \pp^3$ is an Enriques sextic.

\end{thm}
\begin{proof} By   \cite{deg}  the bicanonical image $X$ of $S$ is a surface in
$\pp^{K^2_S}$. Note that necessarily  $p_g(X)=q(X)=0$. Now  assume  that $X$ is
not rational.
 Then   $\deg X\ge 2K_S^2-1$ (see Lemma 1.4 and Remark 1.5 of
\cite{beau}), and, denoting by $d$ the degree of
$\fie$,
 we have $4K^2_S\ge d\deg X\ge d(2K_S^2-1)$, namely $d=2$.

  Denote by $T$ the quotient
 of $S$ by the bicanonical involution. By Theorem 3 and Lemma 7
 of
\cite{xiaocan}, $K^2_S=3$ or $K^2_S= 4$ and there  exists a birational morphism
$T\to \Si$ where $\Si$ is a nodal Enriques surface  with $K^2_S+4$ nodes. The map
$S\to \Si$ is branched on the nodes of $\Si$ and on a  divisor $B$ of $\Si$ with
negligible singularities, contained in the smooth part of $\Si$. Furthermore
$B^2= 2K^2_S$.  The system $|B|$ pulls back to the complete bicanonical system of
$S$, and thus, because
$d=2$, the map determined by $|B|$ is birational and $B$ is nef.

We claim that $|B|$ is base point free. Assume otherwise and suppose first that $|B|$ has a fixed part.
Write
$|B|=|M|+F$, where $M$ and $F$ are the moving part and the fixed part of $B$, respectively.
Since $B$ is contained in the smooth part of $\Si$, also $F$ and the  general
curve in
$|M|$ do not pass through the nodes of $\Si$. The divisor $B$ is nef and
therefore  we have $0\le MB=M^2+MF$, $0\le FB=F^2+MF$ and $MF\geq 1$ (cf. Lemma
2.6 of \cite{ml}). Hence $M^2+MF\le  B^2$, i.e. $M^2\le B^2-MF<B^2=2K^2_S$. This
implies that the image of $\fie$ is a surface of  degree $<2K^2_S$ in
$\pp^{K^2_S}$, which is not possible for an Enriques surface (cf. \cite{codo},
Ch. IV). The same argument shows that, if the general curve in $|B|$ is irreducible then $|B|$ has no base points. Hence $\fie$
is a morphism and $\deg X=2K^2_S$. The surfaces $X$ and $\Si$ have the same minimal
desingularization, which is an Enriques surface.

We now show that  the case $K^2_S=4$ cannot occur. Suppose otherwise. Then the
surface $\Si$ has 8 nodes and thus it is  the surface of Example 1 or the surface
of Example 2. In either case, $\Si$ has  two pencils of elliptic curves $|f_1|$,
$|f_2|$ each containing two double fibres not passing through the nodes of $E$.
Let $2A_1$, respectively $2A_2$,  be such a double fibre of
$f_1$, respectively $f_2$. Recall that $A_1A_2=1$ in  Example 1 and $A_1A_2=2$ in
Example 2. We consider the minimal desingularization
$\eta\colon Y\to \Si$  of $\Si$, and we denote by the same letters the pullbacks
to $Y$ of Cartier divisors of $\Si$.

By Proposition \ref{lattice},
$B\sim_{num}\la_1 A_1 +\la_2A_2$, where $\la_1$, $\la_2$ are integers. The divisor
$B+C_1+\dots +C_8$ is divisible by
$2$ in $\Pic(Y)$  since it is the branch locus of a double cover, and
$C_1+\dots + C_8$ is divisible by 2 by Lemma \ref{even}. Thus $B$ is also
divisible by 2 in $\Pic(Y)$,  hence  $\la_1, \la_2$ are even.
 Finally one has $8=B^2=2\la_1\la_2A_1A_2$. So the only possibility is that
$A_1A_2=1$ and $\la_1=\la_2=2$. Hence
$\Si$ is the surface of Example 1 and $B\sim_{num}f_1+f_2$. Since both $f_1$ and
$f_2$ are $2-$divisible in $\Pic(Y)$, we actually have $B\equiv f_1+f_2$. On the
other hand, the system $|f_1+f_2|$ is not birational  by \cite{codo}, Theorem
4.6.1. Therefore the case $K^2_S=4$ does not occur.
\end{proof}

\begin{rem} In \cite{naie} Daniel Naie constructs surfaces of general type with
$p_g=0$ and $K^2=4$ as double covers of Enriques surfaces. The proof of Theorem
\ref{main} shows that the minimal surfaces $S$ of general type with $p_g=0$ and
$K^2=4$ having an involution $\si$ such that
\begin{itemize}
\item[i)] $S/\si$ is birational to an Enriques surface and
 \item[ii)] the bicanonical map is composed with $\si$ \end{itemize} are
precisely Naie's surfaces.
\end{rem}

\bigskip

\bigskip

\begin{tabbing} 1749-016 Lisboa, PORTUGALxxxxxxxxx\= 56127 Pisa,
ITALY \kill Margarida Mendes Lopes               \> Rita Pardini\\
CMAF \> Dipartimento di Matematica\\
 Universidade de Lisboa \> Universit\a`a di Pisa \\ Av. Prof. Gama Pinto,
2 \> Via
Buonarroti 2\\ 1649-003 Lisboa, PORTUGAL \> 56127 Pisa, ITALY\\
mmlopes@lmc.fc.ul.pt \> pardini@dm.unipi.it
\end{tabbing}

\end{document}